\documentclass[11pt]{amsart}
\usepackage{mathrsfs,amsthm,amssymb,amsmath,url,pst-node,xypic}

\theoremstyle{plain}
    \newtheorem{thm}{Theorem}
    \newtheorem{lem}[thm]{Lemma}
    \newtheorem{prop}[thm]{Proposition}
    \newtheorem{cor}[thm]{Corollary}

    \newtheorem{prob}[thm]{Problem}

\theoremstyle{definition}
    \newtheorem{defn}[thm]{Definition}
    
\theoremstyle{remark}
    \newtheorem{rem}[thm]{Remark}

\newcommand{\nin}{\notin}
\DeclareMathOperator{\ran}{ran} \DeclareMathOperator{\med}{med}
\DeclareMathOperator{\Proj}{Proj}
\renewcommand{\P}{{\mathscr P}}
\newcommand{\inv}{^{-1}}

\newcommand{\un}{^{(n)}}
\newcommand{\uk}{^{(k)}}
\newcommand{\um}{^{(m)}}
\newcommand{\uo}{^{(1)}}
\newcommand{\ut}{^{(2)}}

\newcommand{\cl}[1]{\langle #1 \rangle}

\renewcommand{\O}{{\mathscr O}}
\newcommand{\On}{{\mathscr O}^{(n)}}

\newcommand{\Oo}{{\mathscr O}^{(1)}}
\newcommand{\Ot}{{\mathscr O}^{(2)}}

\DeclareMathOperator{\pol}{Pol}

\newcommand{\C}{{\mathscr C}}
\newcommand{\D}{{\mathscr D}}
\newcommand{\F}{{\mathscr F}}
\newcommand{\J}{{\mathscr J}}

\newcommand{\X}{{\mathscr X}}
\newcommand{\Y}{{\mathscr Y}}
\newcommand{\Z}{{\mathscr Z}}

\newcommand{\I}{{\mathscr I}}
\newcommand{\A}{{\mathscr A}}
\newcommand{\G}{{\mathscr G}}
\renewcommand{\H}{{\mathscr H}}
\newcommand{\M}{{\mathscr M}}
\newcommand{\N}{{\mathscr N}}
\newcommand{\V}{{\mathscr V}}
\renewcommand{\S}{{\mathscr S}}
\newcommand{\Lin}{{\mathscr L}}

\renewcommand{\L}{{\frak L}}
\newcommand{\Pos}{{\frak P}}

\DeclareMathOperator{\spann}{span}

\title[Monoidal intervals]{Monoidal intervals of clones on infinite sets}
\author[M.\,PINSKER]{MICHAEL PINSKER}
\address{Algebra\\TU Wien\\Wiedner Hauptstra\ss e 8-10/104\\A-1040 Wien, Austria}
\email{marula@gmx.at} \urladdr{http://www.dmg.tuwien.ac.at}
\thanks{
Support by the Austrian Science foundation through grant P17627,
and later through the Postdoctoral Fellowship of the Japan Society
for the Promotion of Science (JSPS) is gratefully acknowledged.}

\begin{document}

    \begin{abstract}
        We show that on an infinite set $X$ of cardinality $\kappa$, if $\L$ is the lattice of order ideals
        of some partial order $\Pos$ with smallest element such that $|\Pos|\leq 2^\kappa$, then there is a monoidal interval
        in the clone lattice on $X$
        which is isomorphic to $\L$. In particular, we find that if $\L$ is any chain with smallest element
        which is an algebraic lattice, and if $|\L|\leq 2^\kappa$,
        then $1+\L$ appears as a monoidal interval; also, if $\lambda\leq\kappa$,
        then the power set of $\lambda$ with an additional
        smallest element is a monoidal interval. Concerning
        cardinalities of monoidal intervals these results imply
        that there are monoidal
        intervals of all cardinalities smaller than $2^\kappa$,
        as well as monoidal intervals of cardinality $2^\lambda$,
        for all $\lambda\leq 2^\kappa$.
    \end{abstract}

    \maketitle


    \begin{section}{The problem}
        Let $X$ be a set of cardinality $\kappa$, and
        denote for all $n\geq 1$ the $n$-ary operations on $X$ by $\On$. Then
        $\O=\bigcup_{n\geq 1}\On$ is the set of all finitary operations on $X$.
        A set of operations $\C\subseteq\O$ is called a \emph{clone} iff
        it is closed under composition and contains all
        projections, that is, all functions of the form
        $\pi^n_k(x_1,\ldots,x_n)=x_k$ ($1\leq k\leq n$). The set of
        all clones on $X$ equipped with the order of
        set-theoretical inclusion forms a complete algebraic
        lattice $Cl(X)$.
        After this introductory section, we are going to
        work exclusively with an infinite base set $X$, in which case the cardinality of $Cl(X)$ is
        $2^{2^\kappa}$. For finite $X$ with at least three elements we have
        $|Cl(X)|=2^{\aleph_0}$, and $|Cl(X)|=\aleph_0$ if the base
        set has two elements. Only in the last case the structure
        of the clone lattice has been completely resolved
        \cite{Pos41}. If $X$ has at least three elements, then
        $Cl(X)$ seems to be too large and complicated to be fully
        understood. One approach to this problem is to partition
        the clone lattice into so-called monoidal intervals.

        Let $\G$ be a submonoid of the monoid of unary operations $\Oo$. The set of all
        clones $\C$ with unary part $\G$ (that is, with $\C\uo=\G$, where $\C\uo=\C\cap\Oo$) forms
        an interval $\I_\G$ of the clone lattice; such intervals are referred to as  \emph{monoidal}.
        The smallest element of $\I_\G$ is obviously
        $\cl{\G}$, the clone generated by $\G$ which in this case
        consists of all essentially unary functions (i.e. functions depending on only one variable)
        whose corresponding unary
        function is an element of $\G$. The largest element
        of $\I_\G$ is easily seen to be $\pol(\G)$, defined to contain
        precisely those functions $f\in\O$ for which
        $f(g_1,\ldots,g_{n_f})\in\G$ whenever $g_1,\ldots,g_{n_f}$ are
        functions in $\G$. Functions with this property are called
        \emph{polymorphisms} of $\G$.

        We are interested in the structure of monoidal intervals,
        in particular in the cardinalities monoidal intervals can
        have; this question was first posed by Szendrei \cite{Sze86}.
        One motivation behind this is that if all monoidal
        intervals were similar in some sense, then the problem of
        describing the clone lattice would, up to that similarity, be reduced to the
        description of one monoidal interval, as well as the
        description of the lattice of all submonoids of $\Oo$. If
        on the other hand monoidal intervals could take many forms,
        then this would be another indication that the clone
        lattice is very complicated.

        There is a deeper concept behind the partition of the
        clone lattice into monoidal intervals. If
        $\C,\D\subseteq\O$ are two distinct clones, then there
        exists $n\geq 1$ such that $\C\un\neq\D\un$, where
        $\C\un=\C\cap\On$. Moreover, if this is the case and $m\geq n$, then
        also $\C\um\neq\D\um$. Therefore, we can say that two
        clones are closer the later their $n$-ary parts start to
        differ. More precisely, the function
        $$
            d(\C,\D)=\begin{cases}\frac{1}{2^{n-1}},&\C\neq\D\wedge
            n=\min\{k:\C\uk\neq\D\uk\},\\0,&\C=\D\end{cases}
        $$
        defines a metric on the clone lattice, first introduced by
        Machida \cite{Mac98}. Formulated in this metric, a
        monoidal interval is just an open sphere of radius $1$ in the metric
        space $(Cl(X),d)$. It also makes sense to consider refinements of this
        partition, for example open spheres of radius $\frac{1}{2}$, or
        equivalently sets of clones with identical binary parts;
        they are of the form $[\cl{\H},\pol(\H)]$, where $\H\subseteq\Ot$
        is a set of binary functions closed under composition and containing the two binary projections.

        For a finite base set $X$ it has been observed
        by Rosenberg and Sauer \cite{RSxx} that all intervals are either at
        most countably infinite or
        of size continuum. We shall give a short argument proving this:
        On a finite base set, the clone lattice
        equipped with Machida's metric is homeomorphic to a closed subset of the
        Cantor space $2^\omega$. To see this, notice first that $\O$ is
        countably infinite, and let $(f_i)_{i\in\omega}$ be an
        enumeration of $\O$ with the property that for all $i<j$
        the arity of $f_i$ is not greater than the arity of $f_j$;
        this is possible, since $\On$ is finite for all $n\geq 1$.
        Now we can assign to every set of operations
        $\C\subseteq\O$ a sequence $s(\C)\in 2^\omega$ by defining
        $s(\C)(i)=1$, if $f_i\in\C$, and $s(\C)(i)=0$ otherwise.
        This gives a bijection from the power set $\P(\O)$ of $\O$ onto $2^\omega$, and
        if we extend Machida's metric from the clone lattice to
        $\P(\O)$ (with the same definition), this mapping is
        easily seen to be a homeomorphism. The set of sequences of $2^\omega$ that correspond to
        clones is a
        closed subset of $2^\omega$. Indeed, set for $i\in\omega$
        and $j\in 2$ a set $A_i^j$ to consist of all $s\in 2^\omega$
        with $s(i)=j$; the $A_i^j$ form a clopen subbasis of the topology of
        $2^\omega$. Now the property that $\C\subseteq\O$ contains all
        projections is equivalent to $s(\C)$ being an element of
        $\Lambda_1=\bigcap\{A_i^1: f_i \text{  projection}\}$. Moreover,
        that $\C$ is closed under composition can be stated in the
        language of sequences by saying that $s(\C)$ is an element of
        $$
            \Lambda_2=\bigcap\{(A_{i_0}^0\cup \ldots\cup A_{i_n}^0)\cup A_j^1:
            f_j=f_{i_0}(f_{i_1},\ldots,f_{i_n})\}.
        $$
        Thus $\C\subseteq \O$ is a clone iff $s(\C)$ is an element of
        $\Lambda=\Lambda_1\cap \Lambda_2$, a closed set since both
        $\Lambda_i$ are intersections of closed sets and hence closed
        themselves. Whence, $(Cl(X),d)$ is indeed homeomorphic to a closed subset of $2^\omega$, which immediately
        yields the topological properties of the clone space proven in \cite{Mac98}.\\
        Now if $\C_1\subseteq\C_2\subseteq\O$, then the interval
        $[\C_1,\C_2]$ in the power set of $\O$ corresponds to the
        interval $[s(\C_1),s(\C_2)]$ in $2^\omega$ with the
        pointwise order, a closed set. Therefore it satisfies the
        continuum hypothesis (see \cite{Kec95} for basics of descriptive set theory).
        Also, if $\C_1$ and $\C_2$ are
        clones, then the interval $[\C_1,\C_2]$ in $Cl(X)$
        corresponds to $[s(\C_1),s(\C_2)]\cap\Delta$ in
        $2^\omega$, again a closed set. We conclude that all
        intervals of the clones lattice on a finite set satisfy
        CH. In particular, monoidal intervals can
        only be finite, countably infinite, or of size continuum.

        The same argument does not work for infinite sets, and we
        shall prove that on a countably infinite set
        there exist monoidal intervals of all cardinalities between $\aleph_0$ and $2^{\aleph_0}$.

        Of the possible sizes finite, $\aleph_0$, and
        $2^{\aleph_0}$ for monoidal intervals over a finite set
        with at least three elements, all possibilities occur:
        There must be a monoidal interval of size continuum, since
        there exist only finitely many monoids and
        $|Cl(X)|=2^{\aleph_0}$. Also finite sizes appear, for
        example the interval corresponding to the monoid $\Oo$ is
        of size $|X|+1$ (\cite{Bur67}), and we will see in this paper that the permutation group is
        an example of a monoid whose monoidal interval has only one
        element (for infinite $X$, but the same proof works on finite sets). See
        \cite{PS82}, \cite{Kro95}, \cite{Kro97} for more examples.
        However, for a fixed set, only finitely many finite numbers appear as sizes
        of monoidal intervals, again because there exist only finitely many monoids.
        Krokhin \cite{Kro97} proved that there
        exist countably infinite monoidal intervals over a finite
        set.

        Goldstern and Shelah \cite{GS0x} showed that on a
        countably infinite base set, many monoids define a
        monoidal interval which is as large as the clone lattice
        $(2^{2^{\aleph_0}})$. Starting from this result, we
        investigated the question whether all monoidal intervals
        on infinite sets are that large, and found that the situation is much more diverse.
    \end{section}

    \begin{section}{Results}
        Let $\Pos$ be a partial order. The set of all order ideals (also called lower subsets) on $\Pos$ with the operations
        of set-theoretical intersection and union is a complete algebraic lattice, a sublattice of the
        power set of $\Pos$. We are going to prove the following
        \begin{thm}\label{THM:LIN:powithO}
            Let $X$ be an infinite set of size $\kappa$.
            If $\Pos$ is any partial order with smallest element which has cardinality at most $2^\kappa$,
            and if $\L$ is the lattice of order ideals on $\Pos$, then
            there exists a monoidal interval in the clone lattice over
            $X$ which is
            isomorphic to $\L$.
        \end{thm}

        It is well-known that the class of lattices of order ideals
        is exactly the class of completely distributive algebraic lattices.
        Therefore we have

        \begin{cor}\label{COR:completelyDistr}
            Let $\L$ be a completely distributive algebraic lattice with at
            most $2^\kappa$ completely join irreducible elements. Then there is a monoidal interval in
            $Cl(X)$ isomorphic to $1+\L$, which is to denote $\L$ plus a new smallest element added.
        \end{cor}

        As an immediate consequence we obtain
        \begin{cor}\label{COR:LIN:powerset}
            Let $\lambda\leq 2^\kappa$. Then there is a monoidal
            interval isomorphic to
            $1+\P(\lambda)$, where $1+\P(\lambda)$ is the power set of $\lambda$ with
            a new smallest element added.
        \end{cor}
            Let $\L$ be a chain which is complete as a lattice.
            An element $p\in\L$ is called a \emph{successor} iff there exists $q\in\L$
            with $q<_\L p$ such that the interval $[q,p]_\L$ contains only $p$ and
            $q$. Obviously, the compact elements of $\L$ are exactly the successors
            and the smallest element of $\L$.
            Therefore, $\L$ is a complete algebraic lattice iff the successors are
            unbounded below every $p\in\L$.
        \begin{cor}\label{COR:LIN:chains}
            Let $\L$ be any chain of size at most $2^\kappa$ which
            is a complete algebraic lattice.
            Then there is a monoidal interval isomorphic to
            $1+\L$, which is $\L$ plus a new smallest element added.
        \end{cor}
        \begin{rem}
            Since $Cl(X)$ is an algebraic lattice, all its
            intervals are algebraic. Also, $Cl(X)$ cannot contain any chains larger that $2^\kappa$, since
            there exist only $2^\kappa$ finitary functions on $X$. Hence, these chains are all
            chains which can occur as monoidal intervals (up to the additional smallest element).
        \end{rem}

        \begin{cor}\label{COR:LIN:ordinals}
            If $1\leq \mu\leq 2^\kappa$ is an ordinal, then there is
            a monoidal interval with the order of $\mu$.
        \end{cor}

        \begin{cor}\label{COR:cardinalities}
            On infinite $X$ of size $\kappa$, there exist at least monoidal intervals
            of the following cardinalities:
            \begin{itemize}
                \item{$\lambda$ for all $\lambda\leq 2^\kappa$.}
                \item{$2^\lambda$ for all $\lambda\leq 2^\kappa$.}
            \end{itemize}
        \end{cor}

        Being complete sublattices of the power set of the base set of the partial order,
        the monoidal intervals exposed in our theorem are all
        completely distributive, and therefore still
        quite special lattices. Therefore not surprisingly, they are not all
        monoidal intervals that can appear.

        \begin{prop}
            There exists a nonmodular monoidal interval.
        \end{prop}
        \begin{proof}
            Let $X$ be linearly ordered, and write $\min(x_1,x_2)$ for the
            minimum function, $\med(x_1,x_2,x_3)$ for the median
            function, and $\max(x_1,x_2)$ for the maximum function
            with respect to that linear order. Denote by $\Proj$ the clone of projections. Then
            
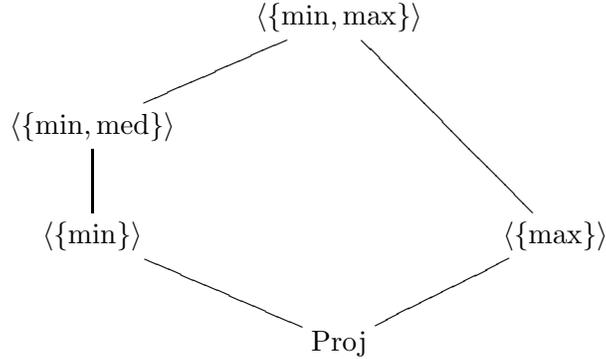
\begin{displaymath}
            \xymatrix{
            &\cl{\{\min,\max\}}\ar@{-}[dl]\ar@{-}[ddr]&\\
            \cl{\{\min,\med\}}\ar@{-}[d]&&\\
            \cl{\{\min\}}\ar@{-}[dr]&&\cl{\{\max\}}\ar@{-}[dl]\\
            &\Proj&
            }
            \end{displaymath}
            is a sublattice of the monoidal interval corresponding
            to the trivial monoid $\{\pi^1_1\}$. That $\cl{\{\min,\med\}}\cap\cl{\{\max\}}=\Proj$
            follows from \cite{Pin032} but is also not difficult to verify.
        \end{proof}

        The fact that monoidal intervals must be algebraic
        lattices with no more than $2^\kappa$ compact (in the clones lattice, this means finitely generated) elements
        is the only restriction for them we know of.
        Therefore we pose the following problem.

        \begin{prob}
            If $\L$ is any algebraic lattice with at most $2^\kappa$ compact elements, is there a
            monoidal interval isomorphic to $\L$?
        \end{prob}

        Concerning cardinalities our theorem leaves the following cases
        open:

        \begin{prob}
            Are the cardinalities of Corollary \ref{COR:cardinalities} all possible sizes of
            monoidal intervals? That is, if $2^\kappa<\lambda<2^{2^\kappa}$ and $\lambda$ is not a cardinality of a
            power set, does there exist a monoidal interval of size $\lambda$?
        \end{prob}

        \subsection{Notation} The smallest clone containing a set
        $\F\subseteq\O$ shall be denoted by $\cl{\F}$; moreover, we write
        $\F^*$ for the set of all functions which arise from
        functions of $\F$ by identification of variables,
        addition of fictitious variables, or permutation
        of variables.
        For $n\geq 1$ we denote the $n$-ary operations on
        $X$ by $\On$; if $\F\subseteq\O$, then $\F\un$ will stand for $\F\cap\On$. We will see $X$
        equipped with a vector space structure; then we write
        $\spann(S)$ for the subspace of $X$ generated by a set of
        vectors $S\subseteq X$. We shall denote the zero vector of $X$ by $0$, and
        use the same symbol for the constant function
        with value $0$. We write $\Lin$ for the set of linear functions on $X$. The sum $f+g$ of two linear functions
        $f,g$ on $X$ is defined pointwise, as is the binary function $f(x)+g(y)$ obtained by the sum
        of two unary functions of different variables. The range of a function $f\in\O$ is given the
        symbol $\ran f$. For a set $Y$ we write
        $\P(Y)$ for the power set of $Y$ and $\P_{fin}(Y)$ for the set of finite subsets of $Y$.
    \end{section}

    \begin{section}{Monoids of linear functions}
        Given any partial order $\Pos$ with $|\Pos|=\lambda\leq 2^\kappa$, we construct a
        monoid $\M$ such that $\I_\M$ is isomorphic to $1+\L$, where $\L$ is the lattice of order
        ideals of $\Pos$.

        Equip $X$  with a vector space structure of dimension $\kappa$
        over any field $K$ of characteristic $\neq 2,3$ and fix a basis $B$ of $X$. Fix moreover
        three distinguished elements $a,b,c \in B$ and write
        $A=B\setminus\{a,b,c\}$.\\ Next we want to introduce a preferably natural notion
        of ``small'' for subsets of $A$; in fact, we are looking for an order ideal
        $\I$
        in $\P(A)$ extending the ideal $\P_{fin}(A)$ which is invariant under permutations of $A$ (i.e., if $S\in\I$ then
        also $\alpha[S]\in\I$ for all permutations $\alpha$ of $X$), such that
        if we factorize $\P(A)$ by this ideal, then the resulting
        partial order has an antichain of length $\lambda$. Since we want to prove our theorem for
        all $\lambda\leq 2^\kappa$, we need the existence of an antichain of length $2^\kappa$, i.e. as large
        as $\P(A)$.
        It is quite obvious that the only order ideals in $\P(A)$
        that are invariant under permutations are the $\I_\xi=\{S\subseteq
        A:|S|<\xi\}$, and the $\J_\xi=\{S\subseteq
        A:|X\setminus S|\geq\xi\}$, where $\xi\leq\kappa$ is a
        cardinal. For $X$ countably infinite, the ideal $\P_{fin}(A)=\I_{\aleph_0}$ satisfies our requirement for
        the antichain.
        For there exists an \emph{almost disjoint} family $\A$ of subsets of
        $A$ of size $2^{\aleph_0}$, meaning that all sets of $\A$
        are infinite and whenever $A_1,A_2\in\A$ are distinct,
        then $A_1\cap A_2$ is finite (see the textbook \cite{Jec02}). The reader interested in
        countably infinite base sets only can imagine this ideal
        in the following. However, it is consistent with ZFC that
        almost disjoint families
        of size $2^\kappa$ fail to exist on uncountable $\kappa$, even if we consider $\I_{\kappa}$ instead
         of $\I_{\aleph_0}$ and replace ``$A_1\cap A_2$ is finite''
         by the weaker
        $``|A_1\cap A_2|<\kappa$''. Moreover, if $\I_{\kappa}$ does not give us an antichain of
        desired length, then the other ideals
        will not work either, so we have to do something
        less elegant: Fix any family $\A\subseteq\P(A)$ of subsets of
        $A$ of cardinality $\kappa$ such that $|\A|=\lambda$, and such that
        $A_1\nsubseteq A_2$ for all distinct $A_1,A_2\in\A$.
        Such families exist; see the textbook \cite[Lemma 7.7]{Jec02} for a proof of
        this. Now we set the ideal $\I$ to consist of all proper
        subsets of sets in $\A$, plus all finite sets, and call the sets of $\I$ \emph{small}. Obviously,
        $\I$ is only an order ideal (no lattice ideal) and quite arbitrary compared
        to the ideal of finite subsets of $A$ which we can use for
        countably infinite $X$. Note also that we had to give up invariance under permutations of $A$;
        however, it will be sufficient
        that if $\alpha$ maps $A_1$ bijectively onto $A_2$, where $A_1,A_2\in\A$,
        and if $S\subseteq A_1$ is small, then
        $\alpha[S]$ is small. Clearly, the sets of $\A$ are not
        elements of $\I$, but their nontrivial intersections are. We index the family $\A$ by the
        elements of $\Pos$: $\A=(A_p)_{p\in\Pos}$.

        The monoid $\M$ we are going to construct will be one of linear functions on the vector space
        $X$, the set of which we denote by $\Lin$. We shall sometimes speak of the
        \emph{support} of a linear
        function $f$, by which we mean the subset of $A$ of those basis vectors which $f$ does not
        send to $0$. The monoid $\M$ will be the union of seven classes of
        functions, plus the zero function. Three classes, namely
        $\N$, $\N'$ and $\N''$, do ``almost nothing'', in the
        sense that they have small support; $\N$ essentially
        guarantees that the polymorphisms $\pol(\M)$ of the monoid $\M$ are sums of
        linear functions, and $\N'$ and $\N''$ are auxiliary functions
        necessary for the monoid to be closed under
        composition. The class $\Phi$ represents the elements of
        the partial order $\Pos$, the class $\Psi$ its order. Finally, the
        classes $\S_\Phi$ and $\S_{\N'}$ ensure that there exist
        nontrivial polymorphisms of the monoid, and that they
        correspond to elements of the partial order.\\

        We start with the set $\N$
        of those linear functions $n\in \Lin$ which
        satisfy the following conditions:
        \begin{itemize}
            \item $n(a)=a$
            \item $n(b)=0$
            \item $n(c)=c$
            \item $n$ has small support.
        \end{itemize}

        Next we add the set $\N'\subseteq\Lin$ consisting of all linear functions $n'$ for which:
        \begin{itemize}
            \item $n'(a)=0$
            \item $n'(b)=0$
            \item $n'(c)=b$
            \item $n'$ has small support
            \item $\ran n'\subseteq\spann(\{b\})$.
        \end{itemize}

        The class $\N''$ contains all $n''\in\Lin$ with
        \begin{itemize}
            \item $n''(a)=a$
            \item $n''(b)=0$
            \item $n''(c)=0$
            \item $n''$ has small support
            \item $\ran n''\subseteq\spann(\{a\})$.
        \end{itemize}

        Observe that all functions $f$ in these three classes have
        small support, and that the range of the functions of
        $\N'$ and $\N''$ is only a one-dimensional subspace of
        $X$.

        Now we define for all $p\in\Pos$ a function
        $\phi_p\in\Lin$ by setting
        \begin{itemize}
            \item $\phi_p(a)=0$
            \item $\phi_p(b)=0$
            \item $\phi_p(c)=b$
            \item $\phi_p(d)=b$ for all $d\in A_p$
            \item $\phi_p(d)=0$ for all other $d\in B$.
        \end{itemize}

        So $\phi_p$ is essentially the characteristic function of
        $A_p$. Observe that
        $\ran\phi_p\subseteq\spann(\{b\})$. We write
        $\Phi=\{\phi_p:p\in\Pos\}$.

        We fix for all $p,q\in\Pos$ with $q\leq_\Pos p$ a function $\psi_{p,q}\in \Lin$ such that
        \begin{itemize}
            \item $\psi_{p,q}$ maps $A_q$ bijectively onto $A_p$
            \item $\psi_{p,q}(a)=a$
            \item $\psi_{p,q}(b)=0$
            \item $\psi_{p,q}(c)=c$
            \item $\psi_{p,q}(d)=0$ for all other $d\in B$
            \item If $q\leq_\Pos r\leq_\Pos p$, then $\psi_{p,r}\circ
            \psi_{r,q}=\psi_{p,q}$.
        \end{itemize}

        This is possible: Let $Y$ be a set of cardinality $\kappa$
        and choose for all $p\in\Pos$ a bijection $\mu_p$ mapping
        $A_p$ onto $Y$. Then setting $\psi_{p,q}(d)=\mu_p\inv\circ
        \mu_q(d)$ for all $d\in A_q$, $\psi_{p,q}(a)=a$, $\psi_{p,q}(c)=c$, and $\psi_{p,q}(d)=0$ for
        all remaining $d\in B$ yields the required functions. We set
        $\Psi=\{\psi_{p,q}:p,q\in\Pos, q\leq_\Pos p\}$. The idea
        behind $\psi_{p,q}$ is that it ``translates''
        the function $\phi_p$ of $\Phi$ into the function
        $\phi_q$, and that such a translation function exists only
        if $q\leq_\Pos p$. More precisely we have

        \begin{lem}\label{LEM:compositionWithPhi}
            Let $\phi_r\in\Phi$ and $\psi_{p,q}\in\Psi$. If $r=p$,
            then $\phi_r\circ\psi_{p,q}=\phi_q$; otherwise,
            $\phi_r\circ\psi_{p,q}\in\N'$.
        \end{lem}
        \begin{proof}
            Assume first that $r=p$. Then
            in the composite $\phi_r\circ\psi_{p,q}$, first $\psi_{p,q}$ maps
            $A_q$ onto $A_p$, and all other vectors of $A$ to $0$,
            and then $\phi_r$ sends $A_r=A_p$ to $b$, so that the composite
            indeed sends $A_q$ to $b$ and all other vectors of $A$ to $0$,
            as does $\phi_q$; one easily checks that also the extra conditions on $a,b,c\in B$
            are satisfied. If on the other hand
            $r\neq p$, then the only basis vectors in $A$ which
            $\phi_r\circ\psi_{p,q}$ does not send to zero are
            those in $\psi_{p,q}\inv [A_r\cap A_p]$, a small set
            since $\psi_{p,q}$ is one-one on its support and by the properties of the family
            $\A$. Moreover,
            $\ran(\phi_r\circ\psi_{p,q})\subseteq\ran\phi_r\subseteq\spann(\{b\})$.
            Hence, since also the respective additional conditions
            on $a,b,c\in B$ are satisfied we have
            $\phi_r\circ\psi_{p,q}\in\N'$.\\
        \end{proof}

        The remaining functions to be added to our monoid are those of the form $\phi_p+n''$,
        where $\phi_p\in\Phi$ and $n''\in\N''$, the set of which we denote by $\S_\Phi$,
        and all functions of the form $n'+n''$, where $n'\in\N'$ and
        $n''\in\N''$; this set we call $\S_{\N'}$. The elements $f$ of
        $\S_\Phi$ and $\S_{\N'}$ both satisfy
        \begin{itemize}
            \item $f(a)=a$
            \item $f(b)=0$
            \item $f(c)=b$.
        \end{itemize}

        We set $\M=\N\cup\N'\cup\N''\cup\Phi\cup\Psi\cup\S_\Phi\cup\S_{\N'}\cup\{0\}$.
        Observe the following properties
        which hold for all $f\in\M$ and which will be useful:

        \begin{itemize}
            \item $f(a)\in\{0,a\}$
            \item $f(b)=0$
            \item $f(c)\in\{0,b,c\}$.
        \end{itemize}

        \begin{lem}\label{LEM:LIN:M_is_a_monoid}
            $\M$ is a monoid.
        \end{lem}
        \begin{proof}
            The following table describes the composition of the
            different classes of functions in $\M$. Here, the meaning of $\X\circ\Y=\Z$ is:
            Whenever $f\in\X$ and
            $g\in\Y$, then $f\circ g\in \Z$.\\
            \begin{center}
            \begin{tabular}{r|c|c|c|c|c|c|c}
                $\circ$      & $\N$      & $\N'$ & $\N''$ & $\Phi$ & $\Psi$                 & $\S_{\Phi}$ & $\S_{\N'}$\\
                \hline
                $\N$         & $\N$      & $0$   & $\N''$ & $0$    & $\N$                   & $\N''$      & $\N''$    \\
                $\N'$        & $\N'$     & $0$   & $0$    & $0$    & $\N'$                  & $0$         & $0$       \\
                $\N''$       & $\N''$    & $0$   & $\N''$ & $0$    & $\N''$                 & $\N''$      & $\N''$    \\
                $\Phi$       & $\N'$     & $0$   & $0$    & $0$    & $\Phi \cup \N'$           & $0$         & $0$       \\
                $\Psi$       & $\N$      & $0$   & $\N''$ & $0$    & $\Psi \cup \N$            & $\N''$      & $\N''$    \\
                $\S_{\Phi}$  & $\S_{\N'}$& $0$   & $\N''$ & $0$    & $\S_{\Phi} \cup\S_{\N'}$& $\N''$      & $\N''$    \\
                $\S_{\N'}$   & $\S_{\N'}$& $0$   & $\N''$ & $0$    & $\S_{\N'}$             & $\N''$      & $\N''$    \\
            \end{tabular}
            \end{center}
            We check the fields of the table. The fact that $\ran
            n'\subseteq\spann(\{b\})$ for all $n'\in \N'$ and $f(b)=0$ for all
            $f\in\M$ yields the $\N'$-column; in the same way we get
            the $\Phi$-column.\\
            If $g=\phi_p+n''\in\S_\Phi$ and $f\in\M$, then $f\circ g=f\circ
            \phi_p+f\circ n''=f\circ n''$, so the $\S_\Phi$-column
            is equal to the $\N''$-column, and the same holds for the
            $\S_{\N'}$-column.\\
            We turn to the $\N$- and $\N''$-columns. The $\S_\Phi$- and the $\S_{\N'}$-row are the sum of the
            $\Phi$- and the $\N'$-row with the $\N''$-row, respectively, since $(f+g)\circ h=(f\circ h)+(g\circ h)$
            for all $f,g,h\in\Oo$. For the other rows of those columns,
            note that if $f,g\in\Lin$ and $g$ has small support, then
            also $f\circ g$ has small support. It is left to the reader to check the conditions on $a,b,c\in B$ and
            on the range for the composites.\\
            It remains to verify the $\Psi$-column. For the first
            row, observe that since all $n\in\N$ have small support and
            since $\psi_{p,q}\inv[S]$ is small for all small $S\subseteq
            A$ and all $\psi_{p,q}\in\Psi$ by the properties of
            $\A$, any composition
            $n\circ\psi_{p,q}$ will have small support. Thus, together
            with the readily checked fact that the extra
            conditions on $a,b,c\in B$ are satisfied we get that
            $n\circ\psi_{p,q}\in\N$. The same argument yields the $\N'$- and
            $\N''$-rows.\\
            The $\Phi$-row is a consequence of Lemma
            \ref{LEM:compositionWithPhi}. Similarly to the proof of that lemma,
            we show that $\psi_{p,s}\circ\psi_{t,q}$ is an element of $\N$
            unless $s=t$, in which case it is $\psi_{p,q}$ by
            construction. Indeed, assume $s\neq t$; then $\psi_{t,q}$
            takes $A_q$ to $A_t$, but $\psi_{p,s}$ has support
            $A_s$; therefore, the composite
            $\psi_{p,s}\circ\psi_{t,q}$ has support $\psi_{t,q}\inv [A_t\cap
            A_s]$, a small set since $\psi_{t,q}$ is injective on its support and
            by the properties of the family $\A$. The conditions on $a,b,c$ for the composite to be in
            $\N$ are left to the reader, and we are done with the $\Psi$-row.\\
            The $\S_\Phi$- and $\S_{\N'}$-rows are the sums of the
            $\N''$-row with the $\Phi$-row and the $\N'$-row
            respectively, by the definitions of $\S_\Phi$ and $\S_{\N'}$.
        \end{proof}

        Recall that if $\F\subseteq\O$, then $\F^*$ consists of all functions
        which arise from functions of $\F$ by
        identification of variables, adding of fictitious variables, as well as by permutation of
        variables. Functions in $\F^*$ are called \emph{polymers}
        of functions on $\F$. Set
        $$
            \V=\{n'(x)+n''(y):n'\in\N', n''\in\N''\}.
        $$
        Moreover,
        define for all $I\subseteq\Pos$ sets of functions
        $$
            \D_I=\{\phi_p(x)+n''(y): p\in I,
            n''\in\N''\}
        $$
        and
        $$
            \C_I=(\M\cup \V\cup\D_I)^*.
        $$
        Observe that $\D_\Pos$ is the set of all functions of the
        form $\phi_p(x)+n''(y)$, where $\phi_p\in\Phi$ and
        $n''\in\N''$.

        \begin{lem}\label{LEM:theCIareClones}
            Let $I\subseteq\Pos$ be an order ideal. Then $\C_I$ is a
            clone in $\I_\M$.
        \end{lem}
        \begin{proof}
            We first show that $\C_I\uo=\M$. Indeed, by its definition the unary functions in $\C_I$ are exactly $\M$
            and those functions which arise when one identifies the two variables of a function in $\V\cup\D_I$.
            If $f\in \V\cup\D_I$, then $f=n'(x)+n''(y)$ or $f=\phi_p(x)+n''(y)$.
            Identifying its variables, we obtain a function of $\S_{\N'}$ in the first and of $\S_{\Phi}$ in
            the second case, and in either case an element of $\M$. Therefore, the unary part of $\C_I$ is
            exactly $\M$ and $\C_I$, if a clone, is indeed an element of $\I_\M$.\\
            $\C_I$ contains $\pi^1_1\in\M$ and therefore all
            projections, as it is by definition closed under the
            addition of fictitious variables.\\ We
            prove that $\C_I$ is closed under composition. To do
            this it suffices to prove that if $f(x_1,\ldots,x_n),
            g(y_1,\ldots,y_m)\in\C_I$, then
            $f(x_1,\ldots,x_{i-1},g(y_1,\ldots,y_m),x_{i+1},\ldots,x_n)\in\C_I$, for all $1\leq i\leq n$.
            Moreover, since $\C_I$ is closed under the addition of fictitious variables,
            we may assume that $f,g$ depend on all of
            their variables, so by the definition of $\C_I$ they are at most
            binary; since within $\C_I$ we can freely permute variables, we can assume $f,g\in\M\cup\V\cup\D_I$.
            Also, since $\C_I$ is by definition closed under identification of variables, we may assume
            that $y_i$ and $x_j$ are different variables, for all $1\leq i\leq m$ and $1\leq j\leq n$.\\
            Let first $f\in\M$. If we substitute any $g\in\M$ for the only variable
            of $f$, then we stay in $\M\subseteq\C_I$ since
             $\M$ is a monoid by Lemma \ref{LEM:LIN:M_is_a_monoid}.
            If $g$ is binary and of the form
            $m'(x)+m''(y)\in\V$, then by Lemma \ref{LEM:LIN:M_is_a_monoid} we have
            $f(m'(x)+m''(y))=f(m'(x))+f(m''(y))=f(m''(y))\in\M^*\subseteq\C_I$,
            since the unary function $f\circ m''\in\M$
            as $\M$ is a monoid.
            Similarly, if $g=\phi_p(x)+m''(y)\in\D_I$ we get $f(\phi_p(x)+m''(y))=
            f(\phi_p(x))+f(m''(y))=f(m''(y))\in\M^*$.\\
            We proceed with the case where $f$ is binary, so $f\in\V\cup\D_I$.
            Assume $f=n'(x)+n''(y)\in\V$, and that we substitute a
            unary $g(z)\in\M$ for $x$. By Lemma \ref{LEM:LIN:M_is_a_monoid}, $n'\circ
            g\in\N'\cup\{0\}$; hence, $f(g(z),y)$ is a function of the form $m'(z)+n''(y)\in\V$
            if $n'\circ g\in\N'$, and the essentially unary function $n''(y)\in\M^*$ if $n'\circ g=0$.
            If we substitute a unary $g(z)\in\M$
            for $y$, then $n''\circ g\in\N''\cup\{0\}$, so
            that again we stay in $\V\cup \M^*$. So say that
            $f=\phi_p(x)+n''(y)\in\D_I$, and that we substitute a
            unary $g(z)\in\M$ for $x$. From Lemma \ref{LEM:LIN:M_is_a_monoid} we know
            that $\phi_p\circ g\in\N'\cup\Phi\cup\{0\}$. If $\phi_p\circ
            g$ vanishes, then we obtain an essentially unary
            function in $(\N'')^*\subseteq\M^*$ for $f(g(z),y)$.
            If
            $\phi_p\circ g\in\N'$, then the sum with
            $n''(y)$ is in $\V$. The interesting case is the one
            where $\phi_p\circ g\in\Phi$; from the
            proof of Lemma \ref{LEM:LIN:M_is_a_monoid} we know
            that this can only happen if $g$ equals some
            $\psi_{s,t}\in\Psi$.
            Moreover, from Lemma \ref{LEM:compositionWithPhi} we
            infer that the composition is only in $\Phi$ if $s=p$,
            and then we have $\phi_p\circ \psi_{p,t}=\phi_t$.
            Hence in this case,
            $f(g(z),y)=\phi_t(z)+n''(y)\in\D_I$ since $t\leq p\in
            I$. To finish the case where we substitute a
            unary function for a variable of a binary function, let $f=\phi_p(x)+n''(y)$
            and substitute $g(z)\in\M$ for $y$. Then, since $n''\circ g\in\N''\cup\{0\}$,
            the result will either be of the form $\phi_p(x)+m''(z)$ and
            thus in $\D_I$, or just $\phi_p(x)\in\M^*$ in case $n''\circ g$ vanishes.\\
            We now substitute binary functions $g(v,w)\in\V\cup\D_I$ into one variable of
            a binary $f(x,y)\in\V\cup\D_I$. Let
            $g(v,w)=m'(v)+m''(w)\in\V$. Since $h\circ m'=0$ for all $h\in\M$, and $f(x,y)$ is of the form
            $f_1(x)+f_2(y)$ for some $f_1,f_2\in\M$,
            and since all involved functions are linear, $m'$ will
            vanish in any substitution with $g$. Therefore
            substituting $g$ is the
            same as substituting only
            an essentially unary function, which we already
            discussed. So let $g(v,w)=\phi_q(v)+m''(w)$. Then again, $h\circ\phi_q=0$ for all $h\in\M$,
            so substitution of $g$ is equivalent to substituting only
            $m''(y)$ and we are done.\\
        \end{proof}

        We now prove that $\cl{\M}$ and the $\C_I$ are the only
        clones in $\I_\M$.

        \begin{lem}
            Let $\G$ be a monoid of linear functions on the vector
            space $X$ which contains the constant function $0$, and let $k\geq 1$ be a natural number.
            If for any finite sequence of vectors $d_1,\ldots,d_k\in X$ there
            exist $e_1,\ldots,e_k\in X$ and
            $h_1,\ldots,h_k\in\G$
            such that $h_j(e_j)=d_j$ and $h_j(e_i)=0$ for all $1\leq i,j\leq k$ with $i\neq j$, then
            all functions in
            $\pol(\G)\uk$ are of the form
            $g_1(x_1)+\ldots+g_k(x_k)$, with $g_1,\ldots,g_k\in\G$.
        \end{lem}
        \begin{proof}
            Let $F(x_1,\ldots,x_k)\in\pol(\G)\uk$. Since $0\in\G$, the
            functions $g_j(x_j)=F(0,\ldots,0,x_j,0,\ldots,0)$ are
            elements of $\G$ for all $1\leq j \leq
            k$. We claim
            $F(d_1,\ldots,d_k)=g_1(d_1)+\ldots+g_k(d_k)$ for all
            $d_1,\ldots,d_k\in X$. Indeed, let $e_1,\ldots,e_k\in
            X$ and $h_1,\ldots,h_k\in\G$ be provided by the assumption of the lemma.
            Then
            $h(x)=F(h_1(x),\ldots,h_k(x))$ is an element of $\G$;
            therefore it is linear. Hence,
            $$
            \begin{aligned}
            h(e_1+\ldots+e_k)&=h(e_1)+\ldots+h(e_k)\\
                             &=F(h_1(e_1),\ldots,h_k(e_1))+\ldots+F(h_1(e_k),\ldots,h_k(e_k))\\
                             &=F(d_1,0,\ldots,0)+\ldots+F(0,\ldots,0,d_k)\\
                             &=g_1(d_1)+\ldots+g_k(d_k).
            \end{aligned}
            $$
            On the other hand,
            $$
            \begin{aligned}
            h(e_1+\ldots+e_k)&=F(h_1(e_1+\ldots+e_k),\ldots,h_k(e_1+\ldots+e_k))\\
                             &=F(h_1(e_1)+\ldots+h_1(e_k),\ldots,h_k(e_1)+\ldots+h_k(e_k))\\
                             &=F(d_1,\ldots,d_k).
            \end{aligned}
            $$
            This proves the lemma.
        \end{proof}

        \begin{lem}\label{LEM:containsNsatisfiesQuasiLinear}
            Let $\G$ be a monoid of linear functions on the vector
            space $X$ which contains $0$. If $\G$ contains $\N$, then the condition
            of the preceding lemma is satisfied for all $k\geq 1$.
        \end{lem}
        \begin{proof}
            Given $d_1,\ldots,d_k\in X$ we choose any distinct
            $e_1,\ldots,e_k\in A$. Now for $1\leq j\leq k$ we define $h_j\in\N$ to
            map $e_j$ to $d_j$, $a$ to $a$, $c$ to $c$, and all remaining basis vectors to
            $0$.
        \end{proof}

        \begin{lem}\label{LEM:sumsOfTwo}
            Let $f,g\in\M$ be nonconstant. If $f+g\in\M$, then
            $f\in\N'\cup \Phi$ and $g\in\N''$ (or the other way
            round).
        \end{lem}
        \begin{proof}
            Observe where the nontrivial functions of $\M$ map $a,c\in
            B$:\\
            \begin{center}
            \begin{tabular}{r|c|c}
                          & $a$ & $c$\\
            \hline
                $\N$      & $a$ & $c$\\
                $\N'$     & $0$ & $b$\\
                $\N''$    & $a$ & $0$\\
                $\Phi$    & $0$ & $b$\\
                $\Psi$    & $a$ & $c$\\
                $\S_\Phi$ & $a$ & $b$\\
                $\S_{\N'}$& $a$ & $b$\\
            \end{tabular}
            \end{center}
            All functions $f\in\M$ satisfy $f(a)\in\{a,0\}$ and
            $f(c)\in\{b,c,0\}$. Hence, if $f+g\in\M$, then
            $f+g(a)=f(a)+g(a)\in\{a,0\}$ and
            $f+g(c)=f(c)+g(c)\in\{b,c,0\}$.
            Since the field $K$ has characteristic $\neq 2,3$ we have that $a+a,b+b,c+c,b+c\nin\{0,a,b,c\}$.
            Thus it can be seen from the table that if
            $f(a)+g(a)\in\{a,0\}$, then at least one of the
            functions must map $a$ to $0$ and thereby be an element of $\N'\cup\Phi$. From
            the condition $f(c)+g(c)\in\{b,c,0\}$ we infer that
            either $f$ or $g$ must map $c$ to $0$ and hence belong to
            $\N''$. This proves the lemma.
        \end{proof}
        \begin{lem}\label{LEM:sumsOfThree}
            Let $f,g,h\in\M$ be nonconstant. Then $f+g+h\nin\M$.
        \end{lem}
        \begin{proof}
            Since $K$ has characteristic $\neq 2,3$ we have that no sum of two or three elements of
            $\{a,b,c\}$ is an element of $\{0,a,b,c\}$. If $f+g+h\in\M$,
            then $f(a)+g(a)+h(a)\in\{a,0\}$. This implies that at
            least two of the three functions have to map
            $a$ to $0$ and therefore belong to $\N'\cup\Phi$.
            Also, $f(c)+g(c)+h(c)\in\{b,c,0\}$, from which we
            conclude that at least two functions must map $c$ to
            $0$ and thus be elements of $\N''$. So one
            function would have to be both in $\N'\cup\Phi$ and in
            $\N''$ which is impossible. Hence, $f+g+h\nin\M$.
        \end{proof}
        \begin{lem}\label{LEM:LIN:PolM_consists_of}
            $\pol(\M)=\C_\Pos$. In particular, all functions in $\pol(\M)$ depend on at most two variables.
        \end{lem}
        \begin{proof}
            Since $\C_\Pos$ is a clone with unary part $\M$ by Lemma
            \ref{LEM:theCIareClones}, we have that
            $\C_\Pos\subseteq\pol(\M)$. To see the other inclusion,
            let $F(x_1,\ldots,x_k)\in\pol(\M)\uk$. Then by Lemma
            \ref{LEM:containsNsatisfiesQuasiLinear},
            $F(x_1,\ldots,x_k)=f_1(x_1)+\ldots+f_k(x_k)$, with $f_i\in\M$, $1\leq i\leq
            k$. We show $F\in\C_\Pos$; since clones are closed under the addition of
            fictitious variables, we may assume that $F$ depends
            on all of its variables, i.e. $f_i$ is nontrivial for
            all $1\leq i\leq k$. If $k=1$, then $F\in\M$, so
            $F\in\C_\Pos$. If $k=2$, then since
            $F(x,x)=(f_1+f_2)(x)$ has to be an element of $\M$, Lemma \ref{LEM:sumsOfTwo} implies that
            up to permutation of variables,
            $F\in\V\cup\D_I\subseteq\C_\Pos$. To conclude, observe that $k\geq
            3$ cannot occur by Lemma \ref{LEM:sumsOfThree}, since
            $F(x,x,x,0,\ldots,0)=f_1(x)+f_2(x)+f_3(x)$ must be an
            element of $\M$ if $F\in\pol(\M)$.
        \end{proof}

        \begin{lem}\label{LEM:containN'+N''impliesN+N''}
            Let $\C$ be a clone containing $\M$ and any function
            of $\V$. Then $\C$ contains $\V$.
        \end{lem}

        \begin{proof}
            Let $n'(x)+n''(y)\in\V\cap\C$, where $n'\in\N'$ and
            $n''\in\N''$, and let $m'(x)+m''(y)$ with $m'\in\N'$ and $m''\in\N''$ be an arbitrary function
            in $\V$. Since $\ran m'=\ran n'=\spann(\{b\})$, there is
            $n_1\in\Lin$ with $m'=n'\circ n_1$. This $n_1$ can be chosen to satisfy $n_1(a)=a$, $n_1(b)=0$, and
            $n_1(c)=c$; also, since $m'$ has small support, we can choose $n_1$ to have small support too.
            Then $n_1\in\N\subseteq\M\subseteq\C$.
            Similarly, there is $n_2\in\N$ such that $m''=n''\circ
            n_2$. Hence, $m'(x)+m''(y)=n'(n_1(x))+n''(n_2(y))\in\C$.
        \end{proof}

        \begin{lem}\label{LEM:containPhi+N''impliesN+N''}
            Let $\C$ be a clone containing $\M$ and any function
            of $\D_\Pos$. Then $\C$ contains $\V$.
        \end{lem}

        \begin{proof}
            Let $\phi_p(x)+n''(y)\in\C\cap\D_\Pos$, where $\phi_p\in\Phi$ and
            $n''\in\N''$.
            Taking any $n\in\N$ we set $n'=\phi_p\circ n\in\N'$. Then $\C$
            contains $n'(x)+n''(y)\in\V$ and hence all
            functions of $\V$ by the preceding lemma.
        \end{proof}

        \begin{lem}\label{LEM:functionsForcedIntoC}
            Let $\C$ be a clone containing $\M$ and a function
            $\phi_p(x)+n''(y)\in\D_\Pos$, where $\phi_p\in\Phi$ and $n''\in\N''$.
            If  $q\leq_\Pos p$ and $m''\in\N''$, then $\C$ contains the function
            $\phi_q(x)+m''(y)$.
        \end{lem}

        \begin{proof}
            As discussed in the proof of Lemma \ref{LEM:containN'+N''impliesN+N''},
            there is $n\in\N$ such that $m''=n''\circ n$. Therefore $\C$ contains
            $\phi_p(\psi_{p,q}(x))+n''(n(y))=\phi_q(x)+m''(y)$.
        \end{proof}

        \begin{prop}\label{PROP:allClonesOfIM}
            If $\C\in\I_\M$ is a clone, then $\C=\M^*=\cl{\M}$, or
            $\C=\C_I$, where $I\subseteq\Pos$ is an order ideal on $\Pos$.
        \end{prop}
        \begin{proof}
            Let $\C\neq \cl{\M}$, that is, $\C$ contains an essentially binary function.
            Set $I=\{p\in\Pos:\exists
            n''\in\N''\,(\phi_p(x)+n''(y)\in\C)\}$. By Lemma
            \ref{LEM:functionsForcedIntoC}, $I$ is an order ideal
            of $\Pos$. We claim $\C=\C_I$. Being elements of $\I_\M$, both $\C$ and
            $\C_I$ have $\M$ as their unary part. Let
            $f(x,y)\in\C\ut$ be essentially binary, i.e. depending on both of its variables;
            then up to permutation of variables,
            $f(x,y)\in\V\cup\D_\Pos$ by Lemma \ref{LEM:LIN:PolM_consists_of}. If $f\in\V$, then
            $f\in\C_I$ by definition of $\C_I$. If $f\in\D_\Pos$, then $f(x,y)=\phi_p(x)+n''(y)$,
            where $p\in\Pos$ and $n''\in\N''$. But then $p\in I$ by
            definition of $I$ and so $f\in\C_I$. Hence,
            $\C\ut\subseteq \C_I\ut$. Because $\C$ contains a
            binary function from $\V\cup\D_\Pos$, Lemmas
            \ref{LEM:containN'+N''impliesN+N''} and
            \ref{LEM:containPhi+N''impliesN+N''} imply
            $\C\ut\supseteq\V$. Also, $\phi_q(x)+m''(y)\in\C\ut$
            for all $q\in I$ and all $m''\in\N''$ by Lemma
            \ref{LEM:functionsForcedIntoC}, so that we have
            $\C\ut\supseteq\C_I\ut$ and thus $\C\ut=\C_I\ut$.
            Lemma \ref{LEM:LIN:PolM_consists_of} implies that
            clones in $\I_\M$ are uniquely determined by their
            binary parts, so that we conclude $\C=\C_I$.
        \end{proof}

        \begin{prop}
            Let $\L$ be the lattice of order ideals on the partial order $\Pos$.
            The monoidal interval $\I_\M$ is isomorphic to
            $1+\L$, which is to denote $\L$ with a new smallest
            element (which corresponds to $\cl{\M}$) added to $\L$.
        \end{prop}
        \begin{proof}
            The mapping $\sigma:1+\L\rightarrow\I_\M$ taking an order
            ideal $I\in\L$ to $\C_I$, as well as the smallest
            element of $1+\L$ to $\cl{\M}$, is obviously a lattice
            homomorphism and injective. By the preceding
            proposition it is also surjective.
        \end{proof}

        \begin{proof}[Proof of Theorem \ref{THM:LIN:powithO}]
            Given a partial order $\Pos$ with smallest element, we
            consider the partial order $\Pos'$ obtained from $\Pos$
            by taking away the smallest element. By the preceding
            proposition, we can construct a monoid $\M$ such that
            $\I_\M$ is isomorphic to $1+\L'$, where $\L'$ is the
            lattice of order ideals on $\Pos'$. Now it is enough to observe that $1+\L'$ is
            isomorphic to the lattice $\L$ of order ideals on $\Pos$.
        \end{proof}

        \begin{proof}[Proof of Corollary \ref{COR:completelyDistr}]
            Let $\L$ be a completely distributive algebraic lattice with at
            most $2^\kappa$ completely join irreducibles. Write
            $\Pos$ for the partial order of completely join
            irreducibles of $\L$ (with the induced order), and write $\L'$ for the lattice of order
            ideals on $\Pos$. The
            mapping
            $$
                \sigma: \quad \begin{matrix} \L &\rightarrow& \L'\\
                p
                &\mapsto&
                \{q\in\Pos:q\leq_\L p\}\end{matrix}
            $$
            is easily seen to be a homomorphism;
            $\sigma$ is bijective because in a completely
            distributive algebraic lattice, every element is a
            join of completely join irreducibles.
        \end{proof}

        \begin{proof}[Proof of Corollary \ref{COR:LIN:powerset}]
            The completely join irreducibles of $\P(\lambda)$ are exactly the singleton
            sets, so there are exactly $\lambda\leq 2^\kappa$ of them and we can
            refer to Corollary \ref{COR:completelyDistr}.
        \end{proof}

        \begin{proof}[Proof of Corollary \ref{COR:LIN:chains}]
            $\L$ is completely distributive algebraic, so
            this is a direct consequence of Corollary \ref{COR:completelyDistr}.
        \end{proof}

        \begin{defn}
            A monoid $\G\subseteq \Oo$ is called \emph{collapsing}
            iff its monoidal interval has only one element, i.e.
            $\cl{\G}=\pol(\G)$.
        \end{defn}

        Denote by $\S$ the monoid of all permutations of $X$.

        \begin{prop}
            $\S$ is collapsing.
        \end{prop}

        \begin{proof}
            Let $f\in\pol(\S)\ut$. Then $\gamma(x)=f(x,x)$ is a permutation.
            Now let $x,y\in X$ be distinct. There exists $z\in X$ with
            $\gamma(z)=f(x,y)$. If $z\nin\{x,y\}$, then we can
            find $\alpha,\beta\in\S$ with $\alpha(x)=x$,
            $\alpha(y)=z$, $\beta(x)=y$, and $\beta(y)=z$. But
            then
            $f(\alpha,\beta)(x)=f(x,y)=f(z,z)=f(\alpha,\beta)(y)$,
            so $f(\alpha,\beta)$ is not a permutation. Thus,
            $z\in\{x,y\}$, and we have shown that
            $f(x,y)\in\{f(x,x),f(y,y)\}$ for all $x,y\in X$.\\
            Next we claim that for all $x,y\in X$, if
            $f(x,y)=f(x,x)$, then $f(y,x)=f(y,y)$. Indeed, consider any permutation $\alpha$
            which has the cycle $(xy)$. Then $f(x,\alpha(x))=f(x,y)=f(x,x)$,
            so $f(y,\alpha(y))=f(y,x)$ has to be different from $f(x,x)$, because otherwise the
            function $\delta(x)=f(x,\alpha(x))\in\S$
            is not
            injective. Hence, $f(y,x)=f(y,y)$.\\
            Assume without loss of generality that $f(a,b)=f(a,a)$ for some distinct $a,b\in X$.
            We claim that $f(a,c)=f(a,a)$ for all $c\in X$.
            For assume not; then $f(a,c)=f(c,c)$ for some $c\in X$, and therefore $f(c,a)=f(a,a)$.
            Let $\beta\in\S$ map $a$ to $b$ and $c$ to $a$.
            Then $f(a,\beta(a))=f(a,b)=f(a,a)$, but also $f(c,\beta(c))=f(c,a)=f(a,a)$, a contradiction since $f$
            preserves $\S$. Hence, $f(a,c)=f(a,a)$ for all $c\in
            X$.\\
            Now if $f(\tilde{a},\tilde{b})\neq f(\tilde{a},\tilde{a})$ for some $\tilde{a},\tilde{b}\in
            X$, then $f(\tilde{a},\tilde{b})=f(\tilde{b},\tilde{b})$ and
            as before we conclude $f(c,\tilde{b})=f(\tilde{b},\tilde{b})$ for
            all $c\in X$. But then
            $f(a,a)=f(a,\tilde{b})=f(\tilde{b},\tilde{b})$, so
            $a=\tilde{b}$; furthermore, $f(a,\tilde{a})=f(\tilde{b},\tilde{a})=
            f(\tilde{a},\tilde{a})\neq f(a,a)$ since we must have $a\neq\tilde{a}$, contradicting
            $f(a,c)=f(a,a)$ for all $c\in X$.\\
            Hence, $f(x,y)=f(x,x)$ for all $x,y\in
            X$ so that $f$ is essentially unary. Therefore, all binary functions of $\pol(\S)$
            are essentially unary. By a result of Grabowski
            \cite{Gra97}, this implies that $\S$ is collapsing. (The
            mentioned result was proved for finite base sets with at least three elements, but
            the same proof works on infinite sets.)
        \end{proof}

        \begin{proof}[Proof of Corollary \ref{COR:LIN:ordinals}]
            The preceding proposition gives us the ordinal $1$.
            For larger ordinals, we can refer to Corollary
            \ref{COR:LIN:chains}.
        \end{proof}
        \begin{proof}[Proof of Corollary \ref{COR:cardinalities}]
            This is the direct consequence of Corollaries
             \ref{COR:LIN:powerset} and \ref{COR:LIN:ordinals}.
        \end{proof}

    \end{section}
    \bibliographystyle{alpha}
    \bibliography{actual_bib}
\end{document}